\documentclass[12pt]{article}
\usepackage{amsfonts,amsmath,amssymb,amsthm,latexsym}

\newtheorem{de}{Definition}[section]
\newtheorem{lm}[de]{Lemma}
\newtheorem{pr}[de]{Proposition}
\newtheorem{co}[de]{Corollary}

\newtheorem{te}[de]{Theorem}

\begin{document}

\title{\textbf{Relative regular modules. Applications to von Neumann regular rings}}
\author{\textbf{L. D\u{a}u\c{s}} \\
{\small Departament of Mathematics, Technical University of Civil Engineering,}\\
{\small Bdul. Lacul Tei 124, RO-020396 Bucharest 2, Romania}\\
{\small e-mail: daus@utcb.ro}
}
\date{}
\maketitle

\begin{abstract}
We use the concept of a regular object with respect to another object in an arbitrary category, defined in \cite{dntd}, in order to obtain the transfer of regularity in the sense of Zelmanowitz between the categories $R-$mod and $S-$mod, when $S$ is an excellent extension of the ring $R$. Consequently, we obtain a result of \cite{ps}: if $S$ is an excellent extension of the ring $R$, then $S$ is von Neumann regular ring if and only if $R$ is also von Neumann regular ring. In the second part, using relative regular modules, we give a new proof of a classical result: the von Neumann regular property of a ring is Morita invariant. \newline
\end{abstract}

\section{Introduction and preliminaries}

${\;\;\;\;\;\;\;}$A ring $R$ is called von Neumann regular if for any $r\in R$ there exists $s\in R$ such that $r=rsr$. Zelmanowitz generalizes this 
concept to modules in \cite{z}: a left module $M$ over the ring $R$ is called regular if for each $m\in M$ there exists $g\in Hom_{R}(M,R)$ such that 
$g(m)m=m$. Since a morphism $f\in Hom_{R}(R,M)$ is uniquely given by an element $m\in M$, one can reformulate the regular module defined by Zelmanowitz as follows: for any $f\in Hom_{R}(R,M)$ there exists $g\in Hom_{R}(M,R)$ such that $f=f\circ g\circ f$. \newline
${\;\;\;\;\;\;\;}$In paper \cite{dntd} it was defined the concept of a regular object with respect to another object in an arbitrary category, which extends the notion of regular module.\newline

\begin{de}
Let $M$ and $U$ be objects of a category $\cal{A}$. We say that $M$ is $U-$regular object if for any morphism $f:U\rightarrow M$ there exists a
morphism $g:M\rightarrow U$ such that $f=f\circ g\circ f$.\newline
\end{de} 

Obviously, if $R$ is a ring and $M$ is a left $R-$module, then $M$ is a regular module if and only if $M$ is $R-$regular in the category $R-$mod.\newline
${\;\;\;\;\;\;\;}$The following result provides a key characterization of regular objects in abelian categories.

\begin{pr}
(\cite[Proposition 3.1]{dntd})Let $M$ and $U$ be objects of an abelian category $\cal{A}$. Then $M$ is $U-$regular if and only if $Ker(f)$ is a direct summand of $U$ and $Im(f)$ is
a direct summand of $M$ for any morphism $f:U\rightarrow M$.\newline
\end{pr}

\section{Excellent extensions}

${\;\;\;\;\;\;\;}$Let $S$ be a ring and let $R$ be a subring of $S$ such that $R$ and $S$ have the same identity 1. We say that the ring $S$ is a \emph{free normalizing extension of }$R$ if there exists a finite set $\{a_{1},a_{2},..,a_{n}\}\subseteq S$ such that $a_{1}=1,\;S=%
\sum_{i=1}^{n}a_{i}R,\;a_{i}R=Ra_{i},$ for all $i=1,..,n$ and $S$ is free with basis $\{a_{1},a_{2},..,a_{n}\}$ as both a left and right 
$R-$module. Further, a free normalizing extension $R\subseteq S$ is called an \emph{excellent extension} if $S$ is left $R-$projective (i.e. if $_{S}N$ is a submodule of $_{S}M$ such that $_{R}N$ is a direct sumand of $_{R}M,$ then $_{S}N$ is a direct sumand of $_{S}M$). \newline
${\;\;\;\;\;\;\;}$The concept of excellent extension was introduced by Passman \cite{p} and named by Bonami \cite{b}. The class of excellent extensions includes $n\times n$ matrix rings and crossed products $R\ast G$ where $G$ is a finite group  with $\left|G\right|^{-1}\in R$. \newline

\begin{lm}
Assume that $M$ is $U-$regular in a category $\cal{A}$. The following assertions hold true:\newline
${\;\;\;\;\;\;\;}$(1) If $U\simeq U'$, then $M$ is $U'-$regular. \newline
${\;\;\;\;\;\;\;}$(2) If $M\simeq M'$, then $M'$ is $U-$regular. 
\end{lm}

The first main result of this paper is:

\begin{te}
Let $S$ be an excellent extension of a ring $R$ and let $M$ be a left $S-$module. Then the following assertions hold:\newline
${\;\;\;\;\;\;\;}$(1) If $_{R}M$ is regular, then $_{S}M$ is regular;\newline
${\;\;\;\;\;\;\;}$(2) If $_{S}M$ is regular and torsion-free module, then $_{R}M$ is regular. 
\end{te}

\begin{proof}
(1) Assume that $_{R}M$ is regular, hence $M$ is $R-$regular in the category $R-$mod. Since $Ra_{i}\simeq R$ as left $R-$modules, by previous lemma we obtain that $M$ is $Ra_{i}-$regular in $R-$mod, for all $i=1,...,n$. According to \cite[Corollary 3.9]{dntd}, it follows that $M$ is 
$\oplus_{i=1}^n Ra_{i}-$regular in $R-$mod, that is $S-$regular in $R-$mod.\newline
${\;\;\;\;\;\;\;}$Now we show, using Proposition 1.2, that $M$ is $S-$regular in $S-$mod. Consider $f:S\rightarrow M$ an arbitrary morphism in $S-$mod.
Then $f$ is a morphism in $R-$mod. Since $M$ is $S-$regular in $R-$mod, then $_{R}Ker(f)$ is a direct summand of $_{R}S$ and $_{R}Im(f)$ is
a direct summand of $_{R}M$. But $S$ is left $R-$projective. It follows that $_{S}Ker(f)$ is a direct summand of $_{S}S$ and $_{S}Im(f)$ is
a direct summand of $_{S}M$, hence $_{S}M$ is regular.\newline
${\;\;\;\;\;\;\;}$(2) We will use the following result obtained by Zelmanowitz: a left $R-$module $M$ is regular if and only if every cyclic submodule of $M$ is projective and a direct summand of $M$ (see \cite[Theorem 2.2]{z}). \newline
${\;\;\;\;\;\;\;}$Let $m$ be an arbitrary element of $M$. Consider the cyclic submodule $Sm$ of $M$. Since $_{S}M$ is a regular module, then $_{S}Sm$ is projective and therefore, by the restriction of scalars, $_{R}Sm$ is projective. On the other hand, since $_{S}M$ is a torsion-free module, by $S=\oplus_{i=1}^n Ra_{i}$, it follows
\begin{equation}
Sm=\oplus_{i=1}^n Ra_{i}m=Rm\oplus(\oplus_{i=2}^n Ra_{i}m).
\end{equation}
Hence $Rm$ is a direct summand of $Sm$ in the category $R-$mod, showing that $_{R}Rm$ is a projective module. On the assumption that $_{S}M$ is regular, we also obtain that $Sm$ is a direct summand of $M$ in the category $S-$mod, hence in the category $R-$mod. Since $Rm$ is a direct summand of $Sm$ in $R-$mod, it follows that $Rm$ is a direct summand of $M$ in the category $R-$mod. Hence $_{R}M$ is regular.\newline
\end{proof}

\begin{co}
Let $S$ be an excellent extension of a ring $R$. Then $R$ is von Neumann regular ring if and only if $S$ is von Neumann regular ring.
\end{co}

\begin{proof}
Suppose that $R$ is von Neumann regular ring, namely $_{R}R$ is a regular module. Since $Ra_{i}\simeq R$ as left $R-$modules, it follows that $_{R}Ra_{i}$ is a regular module, for any $i=1,...,n$. By \cite[Theorem 2.8]{z}, we obtain that $S=\oplus_{i=1}^n Ra_{i}$ is regular as left $R-$module, so, using Theorem 2.2, $S$ is regular as left $S-$module. Hence $S$ is von Neumann regular ring.\newline
${\;\;\;\;\;\;\;}$Conversely, suppose that $S$ is von Neumann regular ring. Let $r$ be an arbitrary element of $R$. Then there exists an element $s\in S$ such that $r=rsr$. Since $S$ is a free extension of the ring $R$ with basis $\{a_{1},a_{2},..,a_{n}\}$, then there exist unique elements $s_{1}, s_{2},..., s_{n}\in R$ such that $s=s_{1}+s_{2}a_{2}+...+s_{n}a_{n}$. Therefore $r=rs_{1}r$, showing that $R$ is von Neumann regular ring. \newline
\end{proof}

\section{Morita equivalence. Morita context}

${\;\;\;\;\;\;\;}$In the following, using again relative regular modules, we give a new proof of a classical result: the von Neumann regular property of a ring is Morita invariant.\newline

\begin{lm}
Let $R$ and $S$ be Morita equivalent rings via the equivalence \newline
$F:R-mod \rightarrow S-mod$ and $M$, $U$ be two left $R-$modules. Then $M$ is 
$U-$regular if and only if $F(M)$ is $F(U)-$regular.
\end{lm} 

\begin{proof}
If $\cal{A}$ and $\cal{B}$ are two categories, it is well-known that a covariant functor $F:\cal{A}\rightarrow\cal{B}$ is an equivalence if and only if it is both faithful and full , and each $B\in\cal{B}$ is isomorphic with some $F(A)$ with $A\in\cal{A}$. Using this result, our lemma is a consequence of
\cite[Proposition 2.4]{dntd}. \newline
\end{proof}

\begin{te}
Let $R$ and $S$ be Morita equivalent rings via inverse equivalences \newline
$F:R-mod \rightarrow S-mod$ and $G:S-mod \rightarrow R-mod$. Then $R$ is von Neumann regular ring if and only if $S$ is von Neumann regular ring.
\end{te}

\begin{proof}
Assume that $S$ is von Neumann regulat ring, hence $S$ is $S-$regular in the category $S-$mod. If we set $Q=G(S)$, by \cite[Theorem 22.1]{af} $_{R}Q$ is a progenerator and $F(Q)\simeq S$. Thus $F(Q)$ is $F(Q)-$regular in $S-$mod. Using Lemma 3.1, it follows that $Q$ is $Q-$regular in $R-$mod. Since $_{R}Q$ is a generator, by \cite[Proposition 17.6]{af} there are a left $R-$module $Q'$ and an integer $n\geq 1$ such that
\begin{equation}
Q^{(n)}\simeq R\oplus Q'
\end{equation}
According to \cite[Theorem 3.8]{dntd} and \cite[Corollary 3.9]{dntd}, $Q^{(n)}$ is $Q^{(n)}-$regular in $R-$mod. Using now \cite[Proposition 2.3]{dntd}, by (2) it follows that $R$ is $R-$regular in $R-$mod, hence $R$ is a von Neumann regular ring.\newline
${\;\;\;\;\;\;\;}$The converse is similar.
\end{proof}

${\;\;\;\;\;\;\;}$A Morita context, denoted by $(R,S,M,N,\phi,\psi)$ consists of two rings $R$, $S$, two bimodules $_{R}M_{S}$, $_{S}N_{R}$ and two bimodules morphisms $\phi:M\otimes_{S}N\rightarrow R$ and $\psi:N\otimes_{R}M\rightarrow S$ such that:
\begin{eqnarray}
\phi(m\otimes n)m'=m\psi(n\otimes m') \\
\psi(n\otimes m)n'=n\phi(m\otimes n')
\end{eqnarray}
for any $m,m'\in M$ and $n,n'\in N$. These conditions ensure that the set $T$ of generalized matrices \newline
\begin{equation*}
\left(
\begin{array}{cc}
r & m \\ 
n & s
\end{array}
\right) , r\in R, s\in S, m\in M, n\in N
\end{equation*}
with ordinary matrix addition and multiplication induced by $\phi$ and $\psi$, forms a ring, called the ring of the Morita context.   \newline
${\;\;\;\;\;\;\;}$A Morita context is strict if both $\phi$ and $\psi$ are epimorphisms. A fundamental theorem of Morita says that the categories of modules over two rings with identity $R$ and $S$ are equivalent if and only if there exists a strict Morita context connecting $R$ and $S$.\newline
${\;\;\;\;\;\;\;}$In order to investigate von Neumann regular rings for Morita contexts we need the following lemma, due to Goodearl. Here we indicate a short proof, using relative regular modules. \newline

\begin{lm}
Let $e_{1},..,e_{n}$ be orthogonal idempotents in a
ring $R$ such that $e_{1}+..+e_{n}=1.$ Then $R$ is regular if and only if
for each $x\in e_{i}Re_{j}$ there exists $y\in e_{j}Re_{i}$ such that $x=xyx$.
\end{lm}

\begin{proof}
Using the isomorphism $Hom_{R}(Re_{i},Re_{j})\simeq e_{i}Re_{j}$ and 
the hypothesis, it follows that $Re_{j}$ is a $Re_{i}-$regular object in $R-mod$, for any $1\leq i,j\leq n$. By \cite[Corollary 3.9]{dntd} we obtain that 
$Re_{j}$ is $Re_{1}\oplus Re_{2}\oplus ... \oplus Re_{n}-$regulat, so $Re_{j}$ is $R-$regular, for any $1\leq j\leq n$. Now, using \cite[Theorem 3.8]{dntd} we get that $R$ is $R-$regular, as desired.  \newline
\end{proof}

\begin{te}
Let $T$ be the ring of a Morita context $(R,S,M,N,\phi,\psi)$ \newline
${\;\;\;\;\;\;\;}$(1) If $T$ is a von Neumann regular ring, then $R$, $S$ are also von Neumann regular rings and $M$, $N$ are regular as left and ringht modules; \newline
${\;\;\;\;\;\;\;}$(2) If the Morita context $(R,S,M,N,\phi,\psi)$ is strict, $R$, $S$ are von Neumann regular rings and $M$, $N$ are regular as left modules, then $T$ is a von Neumann regular ring. 
\end{te}

\begin{proof} (1) Denote 
$e=\left(
\begin{array}{cc}
1 & 0 \\ 
0 & 0
\end{array}
\right)$. 
Then $R\simeq eTe$ and $S\simeq (1-e)T(1-e)$, showing that $R$ and $S$ are von Neumann regular rings. \newline
${\;\;\;\;\;\;\;}$Let $m$ be an arbitrary element in $M$ and consider the matrix 
$\left(
\begin{array}{cc}
0 & m \\ 
0 & 0
\end{array}
\right) \in T$. Then there exists a matrix 
$\left(
\begin{array}{cc}
r' & m' \\ 
n' & s'
\end{array}
\right) \in T$ such that
\begin{equation*}
\left(
\begin{array}{cc}
0 & m \\ 
0 & 0
\end{array}
\right)=
\left(
\begin{array}{cc}
0 & m \\ 
0 & 0
\end{array}
\right)
\left(
\begin{array}{cc}
r' & m' \\ 
n' & s'
\end{array}
\right)
\left(
\begin{array}{cc}
0 & m \\ 
0 & 0
\end{array}
\right)
\end{equation*}
It follows that
\begin{equation*}
\left(
\begin{array}{cc}
0 & m \\ 
0 & 0
\end{array}
\right)=
\left(
\begin{array}{cc}
0 & \phi(m\otimes n')m \\ 
0 & 0
\end{array}
\right)
\end{equation*}
and therefore $m=\phi(m\otimes n')m.$ Define $g:M\rightarrow R$ such that $g(\mu)=\phi(\mu\otimes n')$. Obviously, $g$ is a morphism of left $R-$modules and $g(m)m=m$. Hence $R$ is regular as left $R-$module. On the other hand, since $m=\phi(m\otimes n')m$, by relation (3) it follows that
$m=m\psi(n'\otimes m)$ and thus we obtain that $M$ is regular as right $S-$module.\newline
${\;\;\;\;\;\;\;}$Similarly, we can prove that $N$ is regular as left and right module. \newline
${\;\;\;\;\;\;\;}$(2) We show that $T$ is von Neumann regular ring using previous lemma and the complete set of orthogonal idempotents $\left\{e, 1-e\right\}$, where 
$e=\left(
\begin{array}{cc}
1 & 0 \\ 
0 & 0
\end{array}
\right)$. \newline
${\;\;\;\;\;\;\;}$Let $x\in eT(1-e)$. Then  
$x=\left(
\begin{array}{cc}
0 & m \\ 
0 & 0
\end{array}
\right)$, with $m\in M$. Since $M$ is regular left $R-$module there exists $g\in Hom_{R}(M,R)$ such that $g(m)m=m$. But the morphism $\phi:M\otimes_{S}N\rightarrow R$ is epic, so there exists $n\in N$ such that $\phi(m\otimes n)=g(n)$. Consider
$y=\left(
\begin{array}{cc}
0 & 0 \\ 
n & 0
\end{array}
\right)$. One easly checks that $y\in (1-e)Te$ and  moreover $x=xyx$. \newline 
${\;\;\;\;\;\;\;}$Similarly, if $x\in (1-e)Te$ we can find a matrix $y\in eT(1-e)$ such that $x=xyx$.\newline
${\;\;\;\;\;\;\;}$These facts, together with von Neumann regularity of the rings $R$ and $S$ provide that $T$ is von Neumann regular ring.

\end{proof}
\newpage

\begin{center}
\bf{Acknowledgments}
\end{center}
The author would like to thank Professor Constantin N\u{a}st\u{a}sescu for his important remarks and useful discussions.

\end{document}